# Tree congruence: quantifying similarity between dendrogram topologies

Steven U. Vidovic

*University of Southampton, Library Research Engagement Group, Hartley Library, Highfield Campus, Southampton, UK, SO17 1BJ*

*University of Portsmouth, Palaeobiology Research Group, Burnaby Building, Burnaby Road, Portsmouth, UK, PO1 3QL*

## Abstract

Tree congruence metrics are typically global indices that describe the similarity or dissimilarity between dendrograms. This study principally focuses on topological congruence metrics that quantify similarity between two dendrograms and can give a normalised score between 0 and 1. Specifically, this article describes and tests two metrics the Clade Retention Index (CRI) and the MASTxCF which is derived from the combined information available from a maximum agreement subtree and a strict consensus. The two metrics were developed to study differences between evolutionary trees, but their applications are multidisciplinary and can be used on hierarchical cluster diagrams derived from analyses in science, technology, maths or social sciences disciplines. A comprehensive, but non-exhaustive review of other tree congruence metrics is provided and nine metrics are further analysed. 1,620 pairwise analyses of simulated dendrograms (which could be derived from any type of analysis) were conducted and are compared in Pac-man piechart matrices. Kendall's tau-b is used to demonstrate the concordance of the different metrics and Spearman's rho ranked correlations are used to support these findings. The results support the use of the CRI and MASTxCF as part of a suite of metrics, but it is recommended that permutation metrics such as SPR distances and weighted metrics are disregarded for the specific purpose of measuring similarity.

## Introduction

Dendrograms are branching graphical representations of hierarchical clusters based on similarity. Rooted dendrograms are based on the premise that two given specimens in a sample are more closely related than either one is to a third specimen (i.e. c[a,b]). This study is primarily concerned with rooted dendrograms which require a minimum of 3 specimens to make a statement about their relationships, unrooted dendrograms require a minimum of 4 specimens. Dendrograms are commonly referred to as hierarchical cluster diagrams, trees, phenograms, cladograms and phylogenies. The latter synonyms are due to their use in systematic biology and palaeobiology for depicting evolutionary relationships. Phenograms are constructed from a similarity matrix, while cladograms are chosen heuristically by finding the most parsimonious topology in tree-space or choosing the most sampled topology in tree-space. These graphical representations of the evolutionary tree may be referred to as phylogenies and the terms are often used interchangeably. Dendrograms have become synonymous with evolution but despite that, they have been successfully applied in disciplines including but not limited to linguistics, sociology, scholarly communications (Franceschet 2009), business analytics and manufacturing. There are many reasons a researcher may wish to represent their sample using a dendrogram and they may use different terminology to reflect the sample type or methodology used to estimate the dendrogram. Despite the different applications and terminology used for dendrograms, there are good reasons for all disciplines to compare and describe the topologies consistently. For example, a Systematist may wish to see how adding a new taxon, partitioning data, sequencing different genetic material, or making new morphological observations has affected perceived evolutionary relationships. Biological studies have also compared tree topologies of parasites and hosts (Page 1994) to observe coevolution. Alternatively, manufacturers, developers and distributors may wish to understand their consumers and how they can be concordantly categorised under different key variables to help develop commercial strategies. Furthermore, social scientists could derive a dendrogram from a discriminant function analysis to demonstrate the voting patterns in demographic groups (e.g. Ginsburgh and Noury 2008) and wish to compare how this has changed over time.

In response to this, several methods of measuring or demonstrating the similarity/dissimilarity or congruence/incongruence of dendrograms have been developed. These methods have been reviewed elsewhere in the literature (e.g. Rohlf and Sokal 1981; Rohlf 1982; Swofford 1991; Planet 2006), but, a new method was recently presented (Vidovic 2016; Vidovic and Martill 2017) and due to the scope of the published paper, little was said on the methodology and how it compares to those pre-existing indices. The method termed Clade Retention Index (CRI) was developed to understand topological congruence between evolutionary trees developed by independent research groups who over time borrowed more methods and data from one another. However, here it is recognised that the CRI's potential extends beyond evolutionary studies, therefore, its potential and limitations are explored using simulated dendrograms free from biological references. The CRI is compared to other topological congruence metrics for dendrograms that can be given as indices between 0 and 1. Additionally, a new metric utilising the information available from a maximum agreement subtree and a strict consensus tree derived from the same pair of fundamental dendrograms is proposed and tested.

# Comparison

## Synopsis

This study is primarily concerned with measures of topological congruence between dendrograms that are or can be expressed as an index between 0 and 1, typically with 1 being identical. A commonly used incongruence metric which gives identical trees a score of 0 is also studied (Robinson-Foulds distance). Congruence indices can be broadly divided into three types: permutation; pairwise comparison; and consensus methods. Probabilistic methods also exist, they are discussed below but are not tested due to their lack of descriptive power in a global metric.

Typically, the permutation methods count how many permutations it takes to transform a subject dendrogram into a target dendrogram and derive an index number from that exercise. This can be done across a matrix or using permutation algorithms on the dendrogram itself (e.g. Waterman and Smith 1978). The pairwise comparison methods quantify differences observed between dendrograms (e.g. Robinson and Foulds 1981). And, finally, the consensus methods compute a consensus or compromise tree and measure how resolved it is compared to the most resolved tree possible (e.g. Colless 1980; Mickevich 1978).

Below I give a comprehensive but non-exhaustive account of topological congruence metrics between two fundamental dendrograms.

## Clade Retention Index

The CRI (Vidovic 2016; Vidovic and Martill 2017) is a strict consensus-based metric, which means the indices are the product of the agreement and disagreements observed between two fundamental trees. Fundamental trees lack polytomies, while a strict consensus only preserves relationships shared by each fundamental tree and collapses unshared relationships to the nearest common node, forming polytomies. The CRI method quantifies agreement by counting the shared nodes ($N$) and disagreement by counting the number of individual specimens (or terminal nodes) that terminate in a polytomy ($P$). $P$ is restricted to individual terminal nodes because groups of specimens in a polytomy are far more informative than individual specimens, so this is an intuitive cost. This information is in turn put into context by calculating the maximum possible agreement, which is always the total number of specimens analysed ($T$) less one i.e. $T - 1$. These fundamental components for calculating the level of agreement or disagreement between two tree topologies can be used to formulate equation 1 below.

$$i = -\frac{P-N}{T-1}$$

The maximum agreement score achievable by two fundamental trees is 1, while the maximum disagreement achievable is -1. Therefore, the results of Equation 1 can be normalized to be expressed as an intuitive similarity score given between 0 and 1, where $X$ is 1.

$$CRI = \frac{X}{(Rmax - Rmin) \times (i - Rmin)}$$

$$CRI = \frac{1}{(1 - -1) \times ((-\frac{P-N}{T-1}) - -1)}$$

$$CRI = 0.5 \times ((-\frac{P-N}{T-1}) + 1)$$

## Maximum Agreement Subtree x Consensus Fork

*iCong* (Vienne, Giraud, and Martin 2007) is a metric developed based on the maximum agreement subtree (MAST) which is the largest tree that is agreed upon by two fundamental dendrograms by removing the terminal nodes that are the source of disagreement. In principle, *iCong* is a useful metric because it indicates the specific source of disagreement and quantifies the agreement between two dendrograms. However, in practice, *iCong* can be given as a number greater than 1, which does not help compare between multiple analyses. Here, it is recommended that information available from a MAST is combined with the information available from a strict consensus to indicate the total agreement on both terminal and internal nodes respectively. Like the CRI, $N$ is shared nodes and $T$ is the sample of terminal nodes and here $mT$ is the number of terminal nodes maintained by a maximum agreement subtree. The first part of the equation is effectively Colless' unweighted Consensus Fork (CF), which does not count the root node, so 1 is subtracted from $N$ and 2 is subtracted from $T$. The CF is multiplied by the index obtained by dividing the number of terminal nodes in the MAST by the total number of terminal nodes. The CF is multiplied by the function of the MAST to give a factor of total agreement given the name MASTxCF (*i*) here.

$$i = \frac{N-1}{T-2} \times \frac{mT}{T}$$

## Consensus Fork Index

The Consensus Fork Index (CF) was coined by Colless (1980) (=CI$_c$ Rohlf 1982), to provide a metric given between 0 and 1 of branching points (resolution) on a consensus dendrogram. The CF is calculated by dividing the number of nodes preserved (*N*) in a strict consensus, ignoring the basal branch (*N*-1),  by the maximum possible number of nodes (*T*-2) for a rooted tree.

$$CF = \frac{N-1}{T-2}$$

## Colless' Weighted CF or Mickevich-Platnick Pm

This weighted metric was developed by Colless (1980) concurrently with the CF and it was redeveloped by Mickevich & Platnick (1989) (Swofford 1991). The weighted CF is calculated simply by applying a weight of 1 to each terminal node comprising all subtrees in a consensus at all levels and normalising by the maximum score possible.

$$WCF = \frac{\Sigma N_i}{N_{max}}$$

Colless proposed that the maximum score ($N_{max}$) should be calculated as follows:

$$N_{max} = 0.5(T-1)(T+2)$$

However, tests in PAUP* demonstrate that it is unlikely that Colless' calculation for the maximum score is used.

## Mickevich's Consensus Index

Mickevich's consensus index ($CI_M$) sums the weights given across all subtrees preserved in a consensus tree ($\Sigma N_i$) and normalises it using the divisor $N_{max}$. Where $N_{max}$ is the maximum $\Sigma N_i$ over all possible fully resolved dendrograms derived from a sample of terminal nodes (*T*).

$$CI_M = \frac{\Sigma N_i}{N_{max}}$$

$N_i$ is calculated for each subtree preserved in the consensus at all levels by taking the minimum of two values derived from either taking 1 from the terminal nodes in the subtree or taking the number of terminal nodes in the subtree from the total number of terminal nodes. Therefore, $N_i$ is given as follows:

$$N_i = min\{n_i - 1,\ T - n_i\}$$

For example, a subtree comprising 14 terminal nodes in a tree comprising 21 terminal nodes altogether gives a weight of 7 for that subtree ($N_i$ = *min*{14 – 1, 21 – 14} = 7) and a second subtree of 3 terminal nodes gives a weight of 2 ($N_i$ = *min*{3 – 1, 21 – 3} = 2).

The equation Mickevich gives for resolving $N_{max}$ is as follows:

$$N_{max} = L(\tfrac{T}{2})\, L(\tfrac{T-1}{2})$$

where *L* is the largest integer not larger than its argument.

Continuing the example above, *N*max is given as 100 ($N_{max} = L(10.5) \, L(10) = 10 \times 10 = 100$), therefore, if the sum of the weights given above and all other subtrees is 52, $CI_M$ is 0.520.

Rohlf noted, the maximum $\Sigma N_i$ for $CI_M$ is always achieved by the maximally asymmetrical or unbalanced dendrogram, meaning other topologies that are fully congruent will never score 1. Swofford (1991) demonstrated that a balanced dendrogram can achieve $N_{max}$ but this is only true if said dendrogram is ladderized, thus a balanced bifurcating strict consensus dendrogram demonstrating complete agreement cannot achieve 1.

## Rohlf $CI_1$

Rohlf's $CI_1$ was developed to perform the same role on a consensus tree as $CI_M$, but correcting for tree balance (Rohlf 1982, 137). Rohlf achieved this by adjusting the weights given to subtrees to a simpler monotonic function ($N_i = n_i - 1$) and calculating $N_{max}$ across any possible tree shape by arbitrarily resolving polytomies in a given consensus dendrogram in a consistent manner and calculating the weights across that fully resolved tree (Rohlf 1982, 137–38). This was achieved by adding $\Delta_i$ calculated on every polytomy ($f_i \geq 3$) that requires resolving to $\Sigma N_i$. Where $f_i$ is the number of subtrees and terminal nodes in a given polytomy and $n_i$ are the sizes of those subtrees and terminal nodes given in decreasing order ($n_a, n_b, ..., n_f$).

$$\Delta_i = \sum_{a=2}^{f_i-1} \left[ \sum_{b=1}^{a} n_b - 1 \right]$$

Once the $\Delta_i$ has been calculated for all polytomies, it can be summed and applied in the following equation.

$$CI_1 = \frac{\Sigma N_i}{\Sigma \Delta_i + \Sigma N_i}$$

## Distortion Coefficient

The distortion coefficient was initially introduced by Farris (1973) as the difference between the subject dendrogram and a Matrix Representation with Parsimony (MRP) (Kluge and Farris 1969) of the target dendrogram. Farris calculated this by finding the maximum extra steps on the subject dendrogram for each binary character in the MRP and dividing those by maximum extra steps possible, then averaging the results for all characters. However, here the distortion coefficient is given as it is implemented in TNT, where it appears the ensemble retention index (RI) (Farris 1989) is calculated for a subject dendrogram against an MRP of the target dendrogram. This differs in that the RI is calculated by taking the minimum sum of steps on a particular dendrogram ($S$) from the greatest sum of steps on any dendrogram of the same dimensions ($G$) and dividing that number by the greatest number of steps on any dendrogram less the minimum number of steps on any dendrogram of the same dimensions ($M$).

$$Dist.\ Coef. = RI = \frac{G-S}{G-M}$$

## Robinson-Foulds distance

Robinson-Foulds (R-F) distance (metric) (Robinson and Foulds 1981) might also be referred to as symmetric difference metric, contraction/decontraction metric, or partition metric (Planet 2006). The metric works by sampling two unrooted tree topologies to find groups that appear in one tree but not the other. It is a pairwise sampling method that tests for groups present in tree A, but absent in tree B and present in tree B, but absent in tree A. The first phase of an R-F analysis is to remove all branch relationships present in tree A, but not in tree B, followed by adding branch relationships to tree A, that it lacks, but tree B possesses. The sum of the moves in the two phases is the R-F distance. Because it will cost the same to change tree A into tree B and *vice versa*, no further calculations need to be made. The R-F distance was intended to be expressed as an absolute figure, but it is possible to divide the R-F distance by the maximum possible differences to give an index between 0 and 1. This number needs to be inverted to give 0 as incongruent and 1 as congruent. TNT also calculates the R-F distance as an index number between 0 and 1, where 0 is congruent and 1 is incongruent. Both methods are analysed here.

## SPR or rSPR distances

Subtree pruning and regrafting (SPR) is one of the tree permutation algorithms utilized by cladistic software packages. In TNT, SPR distances are implemented by converting the source fundamental dendrograms to MRP matrices and performing SPR on the subject dendrogram until it matches the target dendrogram (Goloboff 2008). The number of SPR permutations required to change one tree into the other is the same as the converse, therefore, the SPR distance is achieved with a single analysis. Using TNT's algorithms, SPR distance calculations are faster than other tree permutation methods, such as the nearest neighbour interchange (NNI) (Waterman and Smith 1978). However, SPR and NNI distances are limited in that they do not consider the level of disagreement, i.e. two groups with the exact same terminal node composition could be recovered in wildly disparate parts of the dendrogram, but only one permutation (SPR move) would be required to solve this problem. To be able to run two dendrograms from distinct analyses in TNT, they must be reduced to common terminal nodes only.

## Quartets and triplets

Quartets were proposed by Estabrook et al. (1985) specifically to overcome the problem of evolutionary direction. The problem of evolutionary direction is that A+B, +C, +D has the same relationships as C+D, +B, +A if there is no root. The quartets method reduces the tree to sub-groups of four evolutionary units and uses a classification system to compare two trees. The triplets method is similar to the quartet method but is employed for rooted cladograms. The triplet method is most appropriate for this study, but it is obsolete with respect to other methods available and because evolutionary direction is a biologically specific problem, although it is worth noting this method could have applications for hierarchical cluster diagrams derived from a principal component or canonical variate analysis.

## Rand$_c$

Rand (1971) presented one of the earliest measures of similarity between clusters that can be applied to hierarchical clusters and therefore dendrograms. This metric had was rediscovered and adapted several times in the 1970s and 80s (Hubert and Arabie 1985, 194). Rand's *c*, denoted Rand$_c$ here, was originally calculated by scoring all possible clustered pairs with 1 if the pair is together (*T*) in both dendrograms, 1 if the pair is always separate (*S*) in both dendrograms, and 1 if a mix (*m*) of together and separate occurs, else each other category scores 0. The sum of together and separate is then divided by the normalizing factor of the sum of all scores. Rand (1971) gave an example comparing clusters "{(a,b,c),(d,e,f)}" and "{(a,b), (c, d,e),(f)}" as in Table 1.

| score/pair | **ab** | **ac** | **ad** | **ae** | **af** | **bc** | **bd** | **be** | **bf** | **cd** | **ce** | **cf** | **de** | **df** | **ef** | total |
|---|---|---|---|---|---|---|---|---|---|---|---|---|---|---|---|---|
| *T* | 1 | 0 | 0 | 0 | 0 | 0 | 0 | 0 | 0 | 0 | 0 | 0 | 1 | 0 | 0 | 2 |

| | | | | | | | | | | | | | | | |
|---|---|---|---|---|---|---|---|---|---|---|---|---|---|---|---|
| *S* | 0 | 0 | 1 | 1 | 1 | 0 | 1 | 1 | 1 | 0 | 0 | 1 | 0 | 0 | 0 | 7 |
| *m* | 0 | 1 | 0 | 0 | 0 | 1 | 0 | 0 | 0 | 1 | 1 | 0 | 0 | 1 | 1 | 6 |

Table 1. Example given by Rand for scoring ((a,b,c),(d,e,f)) versus ((a,b), (c, d,e),(f)).

Therefore, 9 is divided by 15, giving a score of 0.6. Note, an algebraic function is available for calculating this score (Rand 1971, 847).

Rand (1971) claims that this index scores 1 for absolute similarity and 0 when the clusters have no similarities. However, Rand considers clusters consistently plotting separately a similarity, even if those separations are not consistent themselves e.g. cf in the example given above. Indeed, it is possible for two entirely distinct dendrograms to score greater than 0 (Table 2). For example, ((e,f),((a,b),(c,d))) and ((a,f),((e,c),(b,d))) score 0.533. Therefore, $Rand_c$ does not meet the criteria for inclusion in this study.

| score/pair | ab | ac | ad | ae | af | bc | bd | be | bf | cd | ce | cf | de | df | ef | total |
|---|---|---|---|---|---|---|---|---|---|---|---|---|---|---|---|---|
| *T* | 0 | 0 | 0 | 0 | 0 | 0 | 0 | 0 | 0 | 0 | 0 | 0 | 0 | 0 | 0 | 0 |
| *S* | 0 | 1 | 1 | 1 | 0 | 0 | 0 | 1 | 1 | 0 | 0 | 1 | 1 | 1 | 0 | 8 |
| *m* | 1 | 0 | 0 | 0 | 1 | 1 | 1 | 0 | 0 | 1 | 1 | 0 | 0 | 0 | 1 | 7 |

Table 2. An example of $Rand_c$ scores where there are no shared relationships, but the overall score is greater than 0.

## Clade Concordance Index

Another metric calculated from a strict consensus of two fundamental dendrograms is the Clade Concordance (CC) index (Nixon and Carpenter 1996). However, the CC index is limited in its application as it was specifically developed to measure similarity for evolutionary trees derived from cladistic methods and as such it relies heavily on the data matrix used to estimate the trees. The CC index uses the sum of the greatest character lengths (GL) across each fundamental tree, the length of the shortest (most parsimonious) fundamental tree (PL) and the length of their strict consensus tree (CL) in its calculation. The tree lengths refer to the sum of the character (coded morphological or molecular data in the matrix) transformations across the branches of the trees. As such the CC index can be confounded if the most parsimonious tree is found using ordered states or a weighting procedure i.e. not equal weights for character transformations. Furthermore, the CC index has limited applications even within the field of cladistics because one cannot calculate a meaningful result from fundamental dendrograms derived from distinct matrices. Therefore, in the context of this paper, the CC index is mentioned purely for academic purposes and it will not be further analysed or compared.

$$CC = 1 - \frac{(\sum GLn) - PL}{CL - PL}$$

## Probabilistic methods

Nelson's consensus index ($CI_N$) and the topological incongruence length distance (TILD) offer probabilistic approaches to measuring dendrogram similarity. However, they are analysing how likely it is that those dendrograms would be generated, thus they are not truly topological congruence or incongruence metrics and are not considered further.

# Methodology

In this study, 40 dendrograms are analysed in 180 pairwise comparisons using nine methods. The CF, $CI_M$, W CF, $CI_1$, MASTxCF, R-F, R-F using TNT, CRI and SPR were chosen for comparison. One set of Distortion Coefficient results were also analysed and are available in the underpinning data, but are not presented in the results of this paper because the metric recovers two inconsistent indices depending on which dendrogram is used as the target.

The 40 dendrograms were simulated by creating an MRP of a random dendrogram and retaining suboptimal trees in a heuristic analysis. The 40 dendrograms comprise 4 sets of 10, which each have 11, 21, 31, and 41 terminal nodes respectively. Each set of 10 dendrograms were compared to each other using each topological congruence metric, totalling 45 calculations per set of dendrograms per metric, 180 per metric, and 1,620 overall. The results are presented in Pac-man piechart matrices (Vidovic 2016) for visual comparison, and the results were ranked and compared for concordance/discordance (Kendall's tau-b) and difference (Spearman's rho).

The R scripts, procedures, 40 fundamental trees and their 180 pairwise consensus trees are presented in an independently published dataset.

In practice, to do many of the analyses here, all the terminal node names must be the same and the two dendrogram topologies should be reduced to common terminal nodes only. This has implications at least in systematic biology and palaeobiology where the taxonomic composition of two phylogenies produced by two distinct studies are unlikely to be exactly the same. Helpfully, the taxonomic composition of a dendrogram can be reduced to common terminal nodes using CompPhy (Fiorini et al. 2014).

# Results

All 1,620 analyses, plus the 180 analyses using the Distortion Coefficient are available in the associated dataset. Additionally, the 10 fundamental dendrograms, their consensus trees and the scripts or procedures used to analyse them are also available.
https://doi.org/10.5258/SOTON/D1069

|   |   |   |   |   |   |   |   | Kendall's tau-b |
|---|---|---|---|---|---|---|---|---|
| **CF** | 0.647 | 0.567 | 0.702 | 0.929 | 0.956 | -0.980 | 0.826 | 0.214 |
| 0.826 | **$CI_M$** | 0.837 | 0.748 | 0.627 | 0.668 | -0.634 | 0.466 | 0.333 |
| 0.567 | 0.956 | **W CF** | 0.750 | 0.544 | 0.587 | -0.555 | 0.395 | 0.290 |
| 0.865 | 0.910 | 0.896 | **$CI_1$** | 0.677 | 0.725 | -0.690 | 0.535 | 0.315 |
| 0.987 | 0.823 | 0.729 | 0.859 | **MASTxCF** | 0.913 | -0.911 | 0.792 | 0.292 |
| 0.976 | 0.826 | 0.750 | 0.861 | 0.968 | **R-F** | -0.977 | 0.825 | 0.242 |
| -0.981 | -0.804 | -0.724 | -0.843 | -0.969 | -0.995 | **R-F TNT** | -0.840 | -0.206 |
| 0.931 | 0.635 | 0.548 | 0.698 | 0.920 | 0.933 | -0.947 | **CRI** | 0.158 |
| 0.263 | 0.458 | 0.411 | 0.414 | 0.357 | 0.292 | -0.250 | 0.168 | **SPR** |
| Spearman's rho |   |   |   |   |   |   |   |   |

Table 3. Results of Kendall's tau-b (top right) and Spearman's rho (bottom left) analysis of results (N=1,620) derived from 180 pairwise tree congruence analyses using 9 metrics. The results are presented in a matrix, designed to be read across and down (top right) or down and across (bottom left). CF = Colless' Consensus Fork, $CI_M$ = Mickevich's Consensus Information, W CF = Colless' Weighted Consensus Fork, $CI_1$ = Rohlf's Consensus Information, MASTxCF = Maximum Agreement Subtree x Consensus Fork, R-F = rescaled Robinson-Foulds distance, R-F TNT = Robinson-Foulds distance calculated in TNT, CRI = Clade Retention Index, SPR = Subtree pruning and regrafting similarity index.

## Consensus Fork Index

Colless' CF compares favourably to most congruence metrics available in cladistic packages and presented in the literature. The CF is strongly concordant with both the R-F calculated in TNT (Kendall's tau-b -0.980, Spearman's rho -0.981) and calculated by dividing the number of symmetric differences by the maximum possible number of symmetric differences (Kendall's tau-b 0.956, Spearman's rho 0.976). It also compares favourably with MASTxCF despite the results of the MAST analyses down weighting the CF (Kendall's tau-b 0.929, Spearman's rho 0.987).

## Mickevich's Consensus Information

Mickevich's $CI_M$ is concordant with other congruence metrics, but not strongly associated with them. Results derived from the analysis of simulated trees demonstrate that Mickevich's $CI_M$ is most like other weighted metrics, such as the Weighted CF (Kendall's tau-b 0.837, Spearman's rho 0.956) and Rohlf's $CI_1$ (Kendall's tau-b 0.748, Spearman's rho 0.910). However, other congruence metrics studied exhibit much stronger concordance with each other than Mickevich's $CI_M$ and the Weighted CF do.

## Weighted Consensus Fork

The Weighted Consensus Fork exhibits weaker correlation to other congruence metrics than Mickevich's $CI_M$.

For the Weighted CF, topological incongruence appears to be more strongly weighted against in dendrograms plotting larger samples. In the Weighted CF analyses, those in the sample size of 41 score much lower (e.g. Fig 1) than those in the sample size of 11. By comparison, the CF analysis of sample size 11 exhibits the lowest scores and while the CF analysis of sample size 41 are also low scoring, they are not the lowest. Furthermore, the maximum Weighted CF score (0.425) achieved by two dendrograms with 41 terminal nodes exhibits a greater difference from the median score (0.230) than the same scores for the CF analysis (median = 0.667, max. = 0.795). Whereas, the difference between the maximum and median scores for dendrograms with 11 terminal nodes is much smaller than what is observed in all other metrics with the exception of SPR distances.

## Rohlf's $CI_1$

Rohlf's $CI_1$ is more concordant with all other metrics studied than any of the other weighted congruence metrics. However, Rohlf's $CI_1$ is less well correlated with either Mickevich's $CI_M$ or the Weighted CF than those two are with each other. However, like the Weighted CF and Mickevich's $CI_M$, Rohlf's $CI_1$ appears to down weight observed incongruence more heavily in larger samples, but not to the same extent observed in the former two. Four of the pairwise comparisons of the dendrograms with 11 terminal nodes score less than the lowest score for the sample size of 41 using Rohlf's $CI_1$. This is compared to 19 using CF, but 0 for both Mickevich's CI and the Weighted CF.

## Maximum Agreement Subtree x Consensus Fork

The MASTxCF is among the most concordant with all other tree congruence metrics. Spearman's rho shows MASTxCF to be exceptionally well correlated with the CF, R-F, R-F as calculated by TNT (negatively correlated), CRI, and the Distortion Coefficient (0.987, 0.968, -0.969, 0.920, 0.949 respectively). Kendall's tau-b, which handles ties, finds CF (0.929) and both R-F metrics (0.913 and -0.911) to be more concordant with MASTxCF than the CRI or the Distortion Coefficient (0.792 and 0.841 respectively).

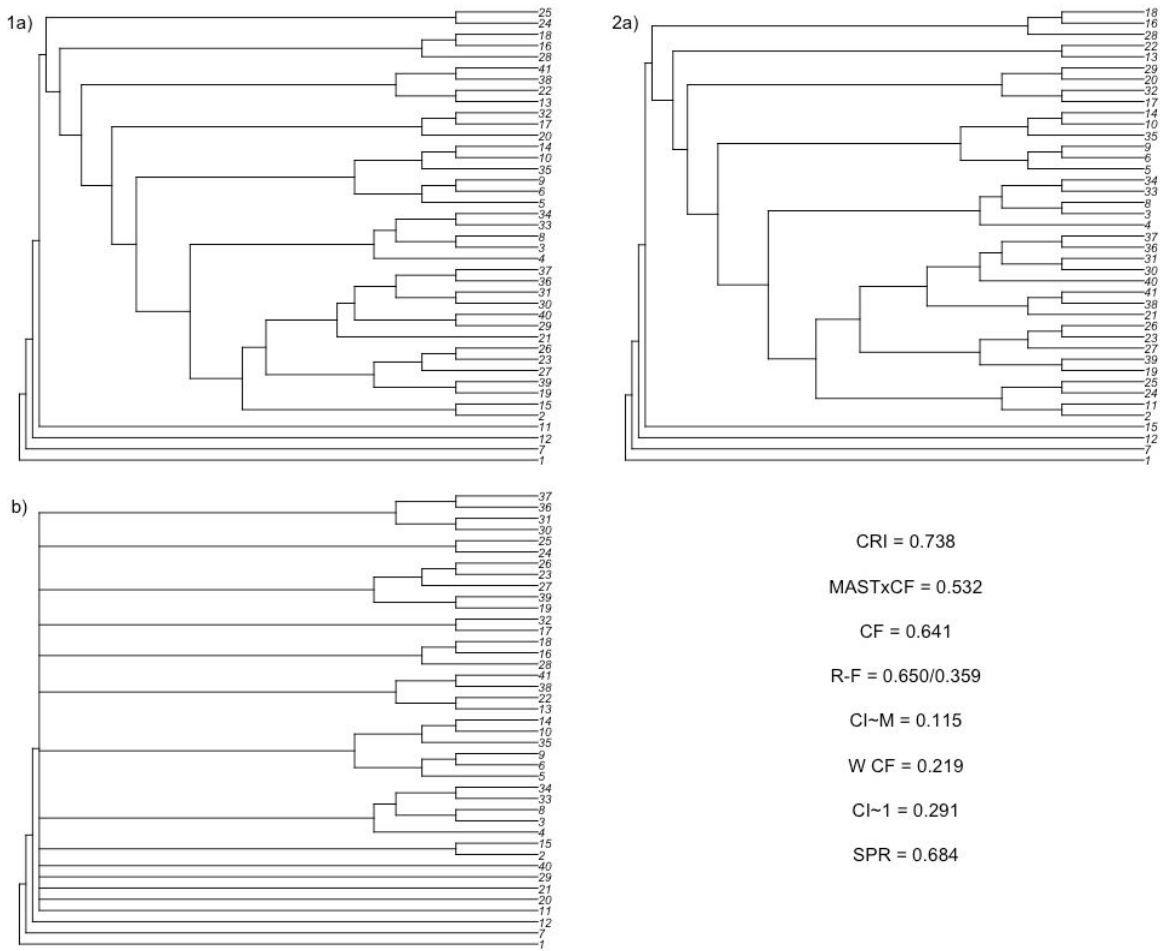

Figure 1. Simulated dendrograms 1a) 3 and 2a) 6 from the sample size of 41, their b) strict consensus and their scores (bottom right).

## Robinson-Foulds Distance

The R-F distances as calculated by TNT and as calculated here find highly similar, but different scores (Kendall's tau-b -0.977, Spearman's rho -0.995). It is unclear how TNT arrives at a different answer, but it is clear that it makes little difference. The R-F distances are among the most concordant with all other tree congruence metrics and even compare favourably with the weighted metrics, with the R-F (as it is calculated here) marginally more concordant with the weighted metrics than MASTxCF (see table 1).

## Clade Retention Index

The CRI is concordant with other metrics but less so than the CF, MASTxCF and R-F. The analyses of dendrograms sampling 11 terminal nodes receive lower scores relative to the larger samples (Fig. 1 & 2. While CF, MASTxCF and R-F all find the sample of 11 terminal nodes to score the lowest out of all the samples, this is more pronounced when using the CRI. Furthermore, when using the CRI, fewer tests of the sample size of 11 score more highly than the sample size of 41 by comparison. However, the sample sizes of 21 and 31 typically score more highly than the sample size of 41 and the highest score for the sample size of 11 is greater than the highest score for 41. Therefore, the smaller samples are not being strongly biased like the larger samples when using a weighted tree congruence metric.

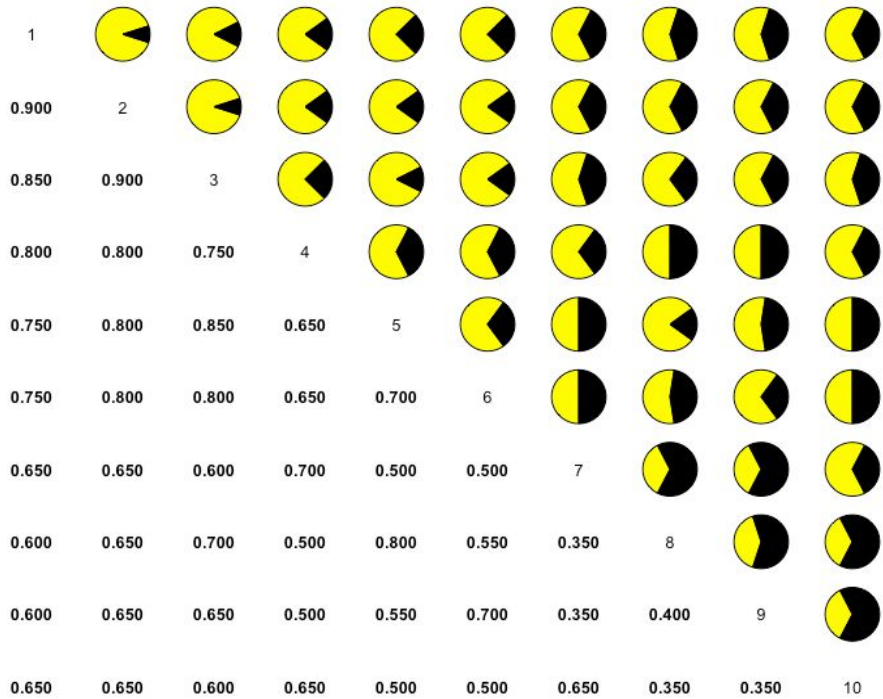

Figure 2. CRI analysis of 10 dendrograms with 11 terminal nodes each. The lower left displays the index numbers derived from 45 consensus trees and the upper right displays the corresponding Pac-man pie charts.

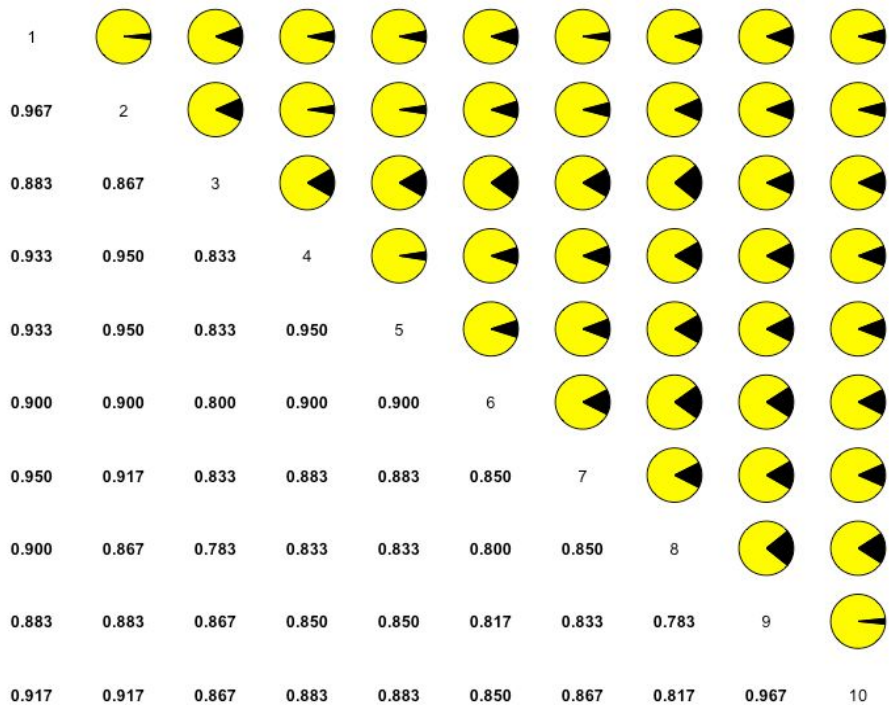

Figure 3. CRI analysis of 10 dendrograms with 31 terminal nodes each. The lower left displays the index numbers derived from 45 consensus trees and the upper right displays the corresponding Pac-man pie charts.

## SPR distances

The SPR distances are by far the least concordant scores compared with all other metrics. The greatest similarity the SPR method has is with Mickevich's CI, but this is weak and largely discordant

(Kendall's tau-b 0.333, Spearman's rho 0.458). Compared with the unweighted metrics, the larger samples score lower relative to the smaller samples. In fact, the sample size of 11 is generally the second-best scoring for SPR distances, while that same sample has among the lowest scores given by other metrics – even the weighted metrics that preferentially down weight larger samples mostly find lower scores for the sample size of 11 (Fig. 1) than those of 21 and 31 (Fig.4).

# Discussion

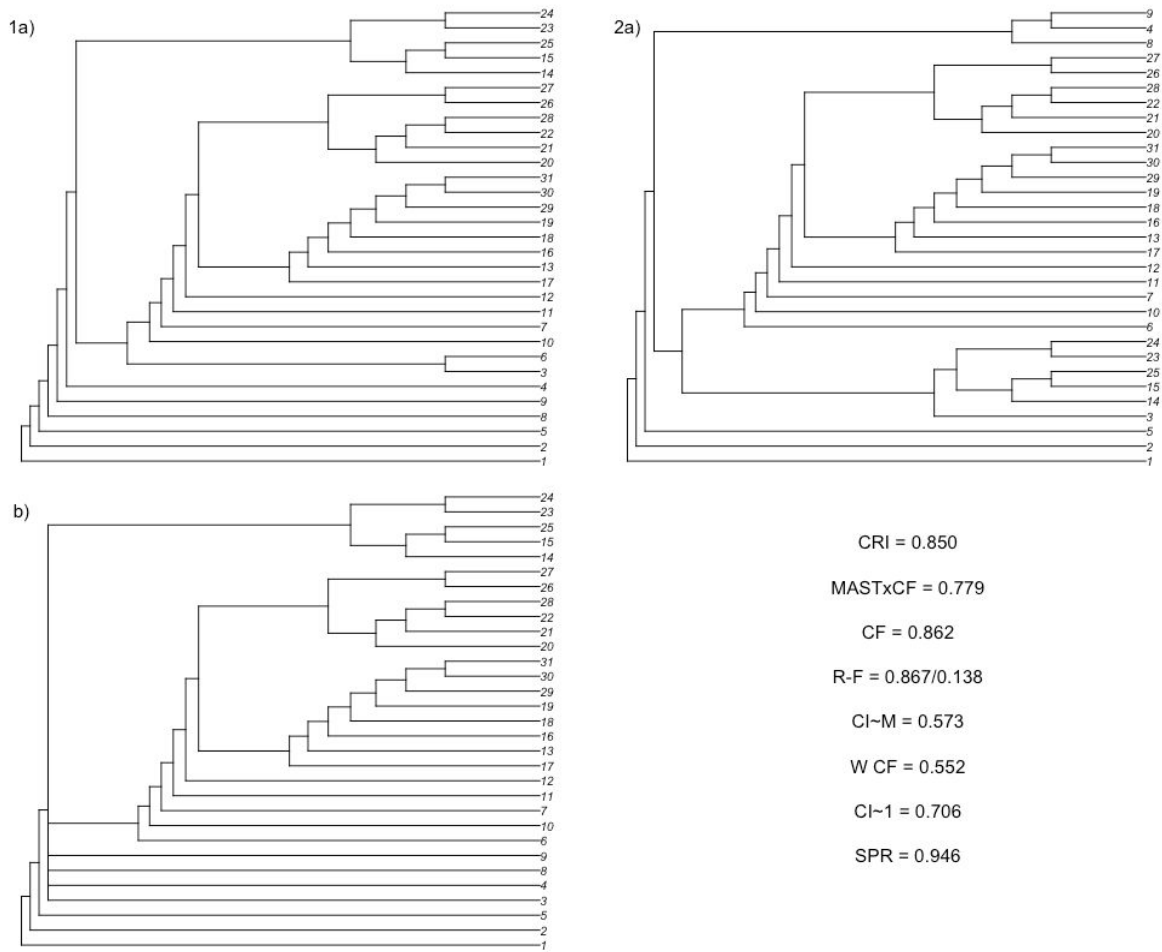

Figure 4. Simulated dendrograms 1a) 5 and 2a) 9 from the sample size of 31, their b) strict consensus and scores (bottom right).

Colless' Consensus Fork index was largely forgotten in the literature and despite being included in a suite of consensus metrics in the cladistics software package PAUP*, the software's developer consigned the CF to "the scrap heap" (Swofford 1991, p. 308). Swofford disregarded the CF because it was unable to differentiate between two example consensus trees with 12 taxa and two nodes each. Swofford (1991) argued that a consensus tree with polytomies in each of its two nested clades preserved more information about relationships than a consensus tree with fewer terminal nodes in a well-resolved pair of nested clades. Here, I agree with that assessment, but also recognise that Swofford's example of a less informative consensus tree demonstrates that there is complete agreement about the relationships of those three resolved terminal nodes, while a far greater number of relationships among the sample could explain the consensus with nested polytomies. In other words, Swofford (1991) was more or equally concerned with preserved information in the consensus tree as he was with topological congruence, leading to misplaced criticism. Here, Kendall's tau-b and Spearman's rho tests indicate that Colless' CF produces results broadly in agreement with commonly preferred metrics, such as the Robinson-Foulds distance. Thus, the CF may have been prematurely dismissed.

Unlike the CF, the CRI and MASTxCF are capable of demonstrating topological congruence and information loss as a result of disagreement. Polytomous terminal nodes and the disagreement between the two fundamental dendrograms that cause them are indicators of incongruence. In response to those tenets, I down weight scores based on preserved nodes by combining information from the composition of polytomies for the CRI and information about how many wildcard terminal nodes are causing disagreement from a MAST for the MASTxCF. Thus, the CRI and MASTxCF may score greater if a few terminal nodes present a well agreed upon relationship, compared to many nested polytomous terminal nodes.

The CRI is concordant with other metrics analysed in this study, but not strongly so. As such, the CRI is telling us something slightly different and used alongside other metrics to interrogate the topological congruence of two fundamental dendrograms it may prove to be a powerful tool. Notably, the CRI results for the smallest sample are skewed towards 0. This may be because less resolution is likely to occur within a polytomy with fewer sampled specimens, as opposed to any active bias. In other words, it is far less likely a dichotomous relationship will occur in a polytomy if there are fewer specimens capable of demonstrating a relationship with one of the polytomous specimens. Also, it stands to reason when working on smaller sample sizes, fewer differences will have a far greater effect on the overall outcome. While a lack of resolution linked to the sample size is undesirable, such an effect is still informative of the lack of information preserved in a consensus which is in contrast to the pronounced sample size bias observed in the weighted suite of metrics. Therefore, the CRI should be used with caution when studying smaller sample sizes but it can be used effectively alongside other metrics.

MASTxCF performs well in these tests and is logical as a metric. As noted above, MASTxCF is capable of demonstrating topological congruence and shared information. This is done by down weighting the score based on common nodes in a strict consensus by a function of the preserved terminal nodes in a MAST. This means MASTxCF is able to demonstrate a difference where incongruence is as a result of significant disruptive disagreement among few specimens, or a little disagreement among many specimens. In the case of the first scenario, the MASTxCF will recover a greater score because a greater proportion of the dendrograms will be in agreement.

Like MASTxCF, weighted metrics attempt to preserve information from the fundamental dendrograms at the same time as demonstrating topological congruence. Unfortunately, the procedures for weighting clusters of relationships in subtrees introduce their own problems. For instance, N*max* in Mickevich's weighted $CI_M$ is only achievable if the tree is maximally bifurcating and maximally asymmetrical. This means the weighted $CI_M$ can only reach 1 if the fundamental dendrograms are fully bifurcating and asymmetrical, therefore, a result less than 1 may not be indicating any incongruence at all. Additionally, all of the weighted metrics tested here were shown to weight scores more heavily against larger samples. This is observed to be a significant effect and will affect researchers' ability to draw comparisons between studies of different sample sizes.

Similarly, SPR distances fail to demonstrate the level of difference between sets of dendrograms being analysed. The results of analyses done here find that one of the worst performing sets of dendrograms for all other metrics is one of the best performing for SPR distances. Even the weighted metrics that preferentially down weight larger dendrograms are more consistent with other metrics

(table 3) than the SPR distances are! Furthermore, because significant disagreements between two tree topologies can be resolved with very few swaps, maybe as little as one, many tests score the same while other congruence metrics find nuanced differences. Therefore, this type of metric demonstrates very little resolution or truth in terms of topological congruence. However, where subtree pruning and regrafting is used to generate dendrograms this method may illustrate how close two fundamental trees are despite appearing drastically different to the eye or other metrics.

The R-F distance is one of the most favoured tree comparison metrics, possibly due to its ability to compare non-binary trees (Robinson and Foulds 1981). The metric has several limitations, including that structurally similar rooted cladograms with one taxon difference can achieve the highest R-F distance possible for a given taxon number (Böcker, Canzar, and Klau 2013). When the actual R-F number is divided by the maximum possible difference between two dendrograms a number between naught and one is produced, but it is not the same as the number produced by TNT. TNT is not open source, thus it is difficult to know what the R-F distance in TNT relates to. Therefore, it is not recommended that the R-F distance is used without any accompanying metrics, such as the CF, CRI or MASTxCF to confirm its findings.

# Conclusion

> "Unfortunately, none of these indices [CF, $CI_1$, $CI_2$ & Pm] provides what we really need: an index that is sensitive to both agreement among and information content of the original trees, but that also allows us to quantify the relative contributions of each of these aspects to lack of resolution in the consensus. We cannot hope to achieve such precision and versatility in a one-dimensional index, but further work in this area may prove fruitful."
> Swofford 1991, p. 310

Similarity is clearly a multidimensional concept (Fowlkes and Mallows 1983), even when comparing pairs of two-dimensional graphs. Indeed, no one metric can be used to inform a Researcher of topological congruence and information content – nor do I think they should. when a metric is calibrated to represent multiple variables, it can become difficult to understand and interpret the effects of states. Among systematists studying evolution, the R-F distance has become popular, but there are potential problems with the R-F distance (Böcker, Canzar, and Klau 2013) and it certainly shouldn't be used alone. The CF and R-F would make an appropriate pairing to help understand the topological congruence between two fundamental dendrograms. Alternatively, the CRI and MASTxCF used together could help a researcher understand the information content and type of incongruence where it is detected as well as how congruent the dendrograms studied are. Combining the CF with those two metrics would deliver yet more information.

Here, it is recommended that weighted metrics are generally avoided unless comparisons between different analyses of different dimensions and tree shapes are not being made – but in that case, it's not clear why you would want to use a metric at all. In the sense of measuring topological congruence, the Distortion Coefficient and SPR distances should be abandoned altogether. The Distortion Coefficient does approximately indicate congruence, but it is not consistently reproducible because the score can vary depending on the reference tree used and the SPR distances are not scalable, nor are they informative.

Beyond what is discussed and analysed here, there are other ways of testing similarity between dendrograms. This includes methods which are available in the R package, Dendextend (Galili 2015). Such as, the Bk (Fowlkes and Mallows 1983) which uses branch lengths to test similarity at different levels of the dendrogram, and Baker's gamma which tests concordance/association between subsets. These methods are not technically topological congruence metrics and therefore were not tested here, but may prove useful being used alongside the methods presented.

**Dataset**

All data supporting this study are openly available from the University of Southampton repository at https://doi.org/10.5258/SOTON/D1069.

How to cite:

Vidovic, S. U. 2019. Dataset for tree congruence: quantifying similarity between dendrogram topologies. ePrints|Soton. https://doi.org/10.5258/SOTON/D1069